\documentclass[a4paper,12pt]{article}

\usepackage[dvips]{graphics}
\usepackage[matrix,arrow,tips,curve]{xy}
\usepackage[english]{babel}
\usepackage{amsmath}
\usepackage{amssymb,amsthm}
\usepackage{enumerate}
\usepackage{psfrag}
  \psfrag{P}{\fontsize{20}{20}$P$}
 \psfrag{Q}{\fontsize{20}{20}$P^*$}

\newtheorem{thm}{Theorem}
\newtheorem*{thm*}{Theorem}

\newtheorem{lemma}[thm]{Lemma}

\theoremstyle{definition}

\theoremstyle{remark}
\newtheorem{remark}[thm]{Remark}
\newtheorem*{example}{Example}

\newcommand{\la}{\longrightarrow}

\newcommand{\ph}{\varphi}
\newcommand{\pr}[1]{\mathbb{P}^{#1}}

\newcommand{\Z}{\mathbb{Z}}
\newcommand{\Q}{\mathbb{Q}}

\newcommand{\ix}{\iota_X}

\newcommand{\Conv}{\operatorname{Conv}}
\newcommand{\Vol}{\operatorname{Vol}}

\newcommand{\Hom}{\operatorname{Hom}}

\begin{document}
\begin{center}
{\bf\large
The number of vertices of a Fano polytope} \\

\bigskip

\textsc{Cinzia Casagrande}

\bigskip

\end{center}

\bigskip

\noindent Let $X$ be a normal, complex, projective variety of dimension $n$. 
Assume that $X$ is Gorenstein and Fano, namely the anticanonical
divisor $-K_X$ of $X$ is Cartier and ample. 
The \emph{pseudo-index} of
$X$ was introduced in \cite{wisn2} as
$$\ix:=\min\{-K_X\cdot C\,|\,C\text{ rational curve in }X\}\in\Z_{>
  0}.$$
By Mori theory we know that $\ix\leq n+1$ when $X$
is smooth, and  $\ix\leq 2n$ in general
(see for instance \cite[Theorems 3.4 and 3.6]{debarreUT}).

 The object of this paper is to give some bounds on the Picard number
$\rho_X$ of $X$, in terms
 of $n$ and $\iota_X$, when $X$ is \emph{toric} and 
\emph{$\Q$-factorial$\,$\footnote{For every Weil
 divisor $D$ there exists $m\in\Z_{>0}$ such that $mD$ is
 Cartier.}}. More precisely,
we will prove the following:
\begin{thm}
\label{uno}
Let $X$ be a $\Q$-factorial, Gorenstein, toric Fano variety of
dimension $n$, Picard number
$\rho_X$ and pseudo-index $\ix$. Then:
\begin{enumerate}[(i)]
\item \label{B}
$\rho_X\leq 2n$, with equality if and only if $n$ is even and $X\cong
  (S_3)^{\frac{n}{2}}$, where $S_3$ is the blow-up of $\pr{2}$ at
  three non collinear points;
\item \label{GM}
 $\rho_X(\ix-1)\leq n$, with equality if and only if 
$X\cong(\pr{\ix-1})^{\rho_X}$.
\end{enumerate}
\end{thm}
Part ($\ref{GM}$)
of Theorem \ref{uno} has been conjectured in \cite{mukai}
for any
smooth Fano variety~$X$, generalizing a conjecture by S.\ Mukai. 
For such $X$,
($\ref{GM}$) 
is known in the cases
$\ix\geq\frac{1}{2}n+1$ \cite{wisn2,cho,occhetta},
$n\leq 4$ \cite{mukai}, $n=5$ \cite{occhettaGM}, and, provided that
$X$ admits an unsplit covering family
of rational curves,
$\ix\geq\frac{1}{3}n+1$ \cite{occhettaGM}.
For $X$ smooth and
toric, ($\ref{GM}$) 
was already known in the cases $n\leq 7$ or
$\ix\geq\frac{1}{3}n+1$ \cite{mukai}.

For a smooth toric Fano $X$,
($\ref{B}$) was conjectured by V.\ V.\ Batyrev 
(see 
\cite[page 337]{ewald2}) and was already
known to hold up to dimension 5
(for $n\leq 4$ thanks to the classifications
  \cite{bat3,wat,bat2,sato}, and for $n=5$ it is \cite[Theorem
  4.2]{fano}). Recently B.\ Nill \cite{nill} has extended this
  conjecture to the $\Q$-factorial Gorenstein case, and has shown
($\ref{B}$) for a certain class of  
$\Q$-factorial, Gorenstein toric Fano
 varieties (see on page \pageref{precise}).

Observe that the bound in ($\ref{B}$) does not hold for non toric Fano
  varieties, already in dimension two. It is remarkable
that in the non toric case, there is
  no known bound for the Picard number of 
  a smooth Fano variety in terms of its dimension
(at least to our knowledge). 
If $S$ is a surface obtained by blowing-up
  $\pr{2}$ at eight general points, for 
 any even $n$ 
 the variety
  $S^{\frac{n}{2}}$ has Picard number
$\frac{9}{2}n$, and one could conjecture this is the maximum for any
dimension  $n$ (see \cite[page 122]{debarre}).

For $X$ smooth, toric, and Fano,
V.\ E.\ Voskresenski{\u{\i}} and A.\ Klyachko have shown that
$\rho_X\leq n^2-n+1$
\cite[Theorem 1]{VK}; O.\ Debarre has
improved this bound in $\rho_X\leq 2+\sqrt{(2n-1)(n^2-1)}$
\cite[Theorem 8]{debarre}.

The common approach to these questions in the toric case, is
to associate a so-called ``reflexive'' 
polytope to a Gorenstein toric Fano variety, and to combine the
techniques coming from
geometry with the ones coming from the theory of polytopes.

We recall some basic notions on reflexive polytopes and toric Fano
varieties; we refer the
reader to \cite{ewald2}, \cite{debarre} and references therein
for more details.
Let $N\cong\Z^n$ be a lattice and let
 $M:=\Hom_{\Z}(N,\Z)$ be the dual lattice.  Set $N_{\Q}:=N\otimes_{\Z}\Q$
  and
$M_{\Q}:=M\otimes_{\Z}\Q$; for $x\in N_{\Q}$ and
$y\in M_{\Q}$ we denote by $\langle x,y\rangle$ the standard pairing.
For any set of points $x_1,\dotsc,x_r\in N_{\Q}$, we denote by
$\Conv(x_1,\dotsc,x_r)\subset N_{\Q}$ their convex hull.

Let $P\subset N_{\Q}$ be a lattice
polytope of dimension $n$
containing the origin in its interior. 
We denote by $V(P)$ the set of vertices of $P$.
The \emph{dual
polytope} of $P$ is defined as
$$ P^*:=\{y\in M_{\Q}\,|\,\langle x,y\rangle \geq -1\text{ for all
}x\in P\}.$$
$P$ is called a \emph{reflexive polytope} if $P^*$ is a lattice
polytope; if so, also $P^*$ is reflexive and $(P^*)^*=P$. Reflexive
polytopes were introduced in
\cite{batdual};
 their isomorphism
classes 
are in bijection with isomorphism classes of Gorenstein
toric Fano varieties.

Let $P$ be a reflexive polytope of dimension $n$;
we denote by $X_P$ the associated 
$n$-dimensional Gorenstein toric Fano variety. The fan of
$X_P$ is given by the cones over the faces of $P$ in $N_{\Q}$.

Many geometric properties of $X_P$ can be read from $P$. In particular, 
$X_P$ is
 $\Q$-factorial if and only if  $P$ is simplicial\footnote{A polytope 
is simplicial  if the vertices of every
 facet are linearly independent, where a  facet is 
  a proper face of maximal dimension.}, while 
 $X_P$ is smooth  if and only if the vertices of every
 facet
of $P$ are a basis of
the lattice. In this last case, we
 say that $P$ is smooth. 
\begin{example}
\label{ex}
In dimension two,
let $e_1,e_2$ be a basis of $N$, and  $e_1^*,e_2^*$ the
dual basis of $M$. 
Consider
$P:=\Conv(e_1,e_2,-e_1-e_2)$; its dual polytope is
$P^*=\Conv(2e_1^*-e_2^*,
2e^*_2-e^*_1,-e_1^*-e_2^*)$. 
\begin{center}
{}{}\scalebox{0.30}{\includegraphics{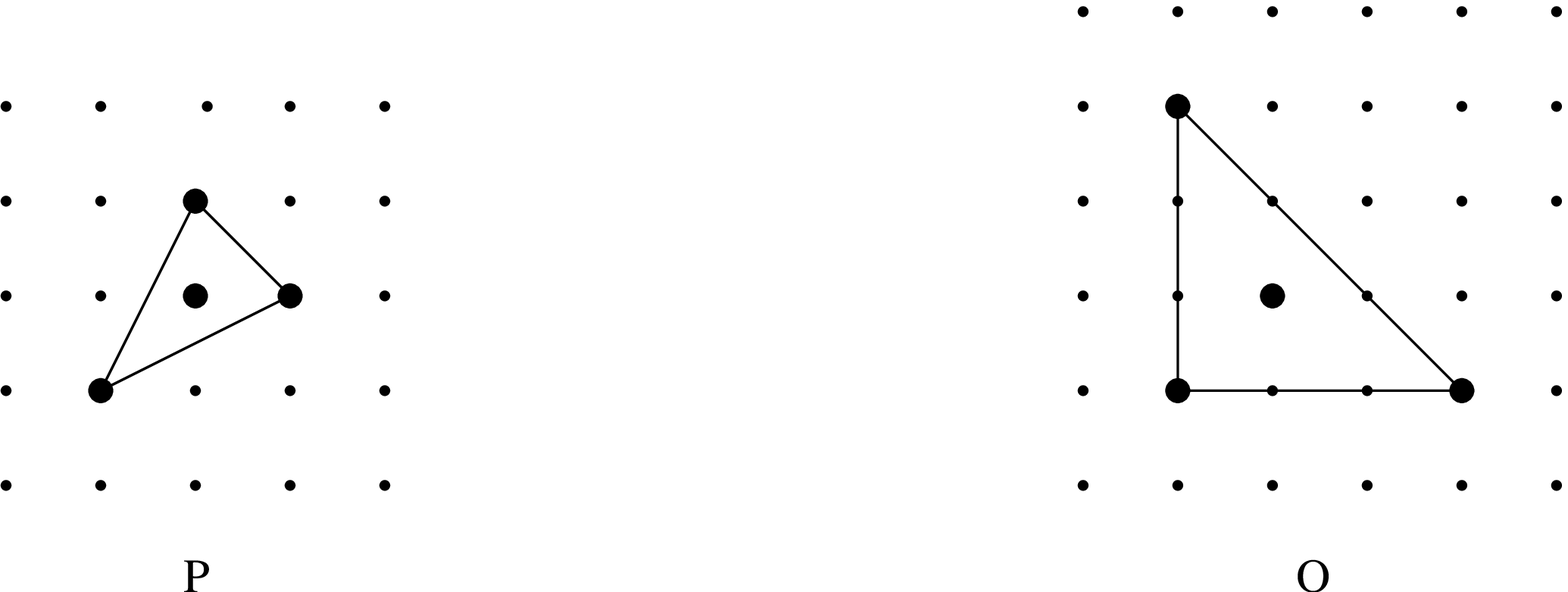}}
\end{center}
Both polytopes are reflexive, the
associated surfaces are  $X_P=\pr{2}$ and $X_{P^*}=\{xyz=w^3\}\subset\pr{3}$.
The surface $X_{P^*}$ has three
singular points of type $A_2$. We have $\iota_{X_P}=3$, while
$\iota_{X_{P^*}}=1$. 
\end{example}
We denote by $|A|$ the cardinality of a finite set $A$. When $P$ is simplicial,
the number of vertices $|V(P)|$ of $P$ is equal to $\rho_{X_P}+n$. 

Recall that there is a bijection between the vertices of $P^*$ and the
facets of $P$; if
$u\in V(P^*)$, we denote by $F_u$ the corresponding
facet of $P$, namely $F_u:=\{x\in P\,|\,\langle x,u\rangle=-1\}$.
We define
$$
\delta_P:=\min\{\,\langle v,u\rangle \,|\, v\in V(P),\,u\in
V(P^*),\,v\not\in F_u\}\,\in\Z_{\geq 0}.$$
The pseudoindex $\iota_{X_P}$ is related to $\delta_P$ as follows. 
\begin{lemma}
\label{pseudoindex}
Let $P$ be a simplicial reflexive polytope and $X_P$ the associated Fano
variety. Then $\iota_{X_P}\leq\delta_P+1$.
If moreover $P$ is smooth, then $\iota_{X_P}=\delta_P+1$.
\end{lemma}
Theorem \ref{uno} is then a consequence of Lemma \ref{pseudoindex} and
of the following:
\begin{thm}
\label{due}
Let $P$ be a simplicial reflexive polytope of dimension $n$. Then:
\begin{enumerate}[(i)]
\item 
\label{duebat}
$|V(P)|\leq 3n$, with equality if and only if $n$ is even and $X_P\cong
  (S_3)^{\frac{n}{2}}$;
\item \label{dueGM} 
 if $\delta_P>0$, then $|V(P)|\leq n+\frac{n}{\,\delta_P}$, 
with equality 
if and only if the following two conditions hold:
\begin{enumerate}[a.]
\item
$P=\Conv(Q_1,\dotsc,Q_r)$ with $r=\frac{n}{\delta_P}$ and each 
$Q_j\subset N_{\Q}$ a reflexive
lattice simplex
of dimension  $\delta_P$, the sum of whose vertices is zero;
\item
if $H_j$ is the linear span of $Q_j$ in $N_{\Q}$, we have 
$N_{\Q}=H_1\oplus\cdots\oplus H_r$.
\end{enumerate}
\end{enumerate}
\end{thm}
In \cite{nill}, Theorem \ref{due} ($\ref{B}$) \label{precise}
is proven for a simplicial reflexive polytope $P$
for which there exists a
vertex $u$ of $P^*$ such that $-u\in P^*$ \cite[Theorem 5.8]{nill}. 
We refer the reader to \cite{nill} for 
a discussion on the number of vertices of a reflexive polytope in the
non simplicial case.
\begin{example}
The polytopes $P$ and $P^*$ of the previous example have
$\delta_P=\delta_{P^*}=2$, and both satisfy
equality in ($\ref{dueGM}$). Notice that $\iota_{X_{P^*}}<\delta_{P^*}+1$.
\end{example}
We will first prove Theorem \ref{due}, then
Lemma~\ref{pseudoindex} and Theorem \ref{uno}.
For the proof, we will use the same technique as
\cite{VK} and 
\cite{debarre}, and also some results from \cite{nill}.

First of all, we need a property of pairs of vertices $v\in V(P)$ and
$u\in V(P^*)$ with $\langle v,u\rangle=\delta_P$.

Let $P$ be a simplicial polytope of dimension $n$.
We say that a vertex $v$ is \emph{adjacent} to a facet
$F=\Conv(e_1,\dotsc,e_n)$ if
$\Conv(v,e_1\dotsc,\check{e}_i,\dotsc,e_n)$ is a facet of $P$ for some
$i=1,\dotsc,n$. 
\begin{lemma} 
\label{zero}
Let $P$ be a simplicial reflexive polytope and $v\in V(P)$, $u\in
V(P^*)$ such that $\langle v,u\rangle=\delta_P$.
Then $v$ is adjacent to $F_{u}$. 
\end{lemma}
\begin{proof}
This property is shown in \cite[Remark 5(2)]{debarre} 
and \cite[Lemma 5.5]{nill} in the case $\delta_P=0$. 
The same proof works for the general case.
\end{proof}
\begin{proof}[Proof of Theorem \ref{due}]
First of all, observe that for any
$u\in V(P^*)$ we have
\begin{equation}\label{bounds}
\left|\left\{v\in V(P)\,|\,\langle v,u\rangle=-1
\right\}\right|=n\quad\text{and}\quad\left|\left\{v\in
  V(P)\,|\,\langle v,u\rangle=0 
\right\}\right|\leq n.\end{equation}
In fact, since $P$ is simplicial, the facet $F_u$ contains
$n$ vertices.
Moreover, if $\langle v,u\rangle=0$, then $\delta_P=0$, and
by Lemma \ref{zero} we know that $v$ is adjacent to $F_{u}$.
Again, since $P$ is simplicial, $F_u$ has at most $n$ adjacent
 vertices, and we get \eqref{bounds}.

The origin lies in the interior of $P^*$, so we can write a
relation
\begin{equation}
\label{sum}
m_1u_1+\cdots+m_hu_h=0
\end{equation}
where $h>0$, $u_1,\dotsc,u_h$ are vertices of $P^*$, and
$m_1,\dotsc,m_h$ are positive integers.
Set $I:=\{1,\dotsc,h\}$ and $M:=\sum_{i\in I}m_i$.
For any vertex $v$ of $P$ define 
$$ A(v):=\{i\in I\,|\,\langle v,u_i\rangle=-1\}\quad\text{and} \quad
B(v):=\{i\in I\,|\,\langle v,u_i\rangle = 0\}.$$
Then observe that $\langle v,u_i\rangle\geq 1$ for any $i\not\in
A(v)\cup B(v)$.
So for every $v\in V(P)$ we have
\begin{align*}
0&=\sum_{i\in I}
m_i\langle v,u_i\rangle 
=-\sum_{i\in A(v)}m_i+\sum_{i\not\in A(v)\cup B(v)} 
m_i\langle v,u_i\rangle \\&\geq -\sum_{i\in A(v)}m_i+\sum_{i\not\in
  A(v)\cup B(v)} m_i= 
M-2\sum_{i\in A(v)}m_i-\sum_{i\in B(v)}m_i.\end{align*} 
Summing over all vertices of $P$ we get
\begin{align*}
M|V(P)| &\leq 2\sum_{v\in V(P)}
\sum_{i\in A(v)}m_i
+\sum_{v\in V(P)}\sum_{i\in B(v)}m_i\\
=&2\sum_{i\in I}m_i\left|\left\{v\in V(P)\,|\,\langle v,u_i\rangle=-1
\right\}\right|+
\sum_{i\in I}m_i|\{v\in V(P)\,|\,\langle v,u_i\rangle=0\}|
\end{align*}
and using \eqref{bounds} this gives $|V(P)|\leq 3n$.

Assume that $|V(P)|=3n$. Then all inequalities above are equalities;
in particular, for any $v$ and $u_i$ such that $\langle
v,u_i\rangle >0$, we must have $\langle
v,u_i\rangle =1$. 
Observe now that we can choose a relation
as \eqref{sum}
involving all vertices of $P^*$, namely with $h=|V(P^*)|$ (see Remark
\ref{relations}), so $\langle v,u\rangle\in\{-1,0,1\}$ 
for every $v\in V(P)$ and $u\in V(P^*)$. Then
$P$ and  $P^*$ are centrally symmetric.

Smooth centrally symmetric reflexive polytopes are classified
in \cite[Theorem 6]{VK}, and the only case with $3n$ vertices is for $n$ 
even and  $X_P\cong
  (S_3)^{\frac{n}{2}}$.
 For the general case, we apply
\cite[Theorem 5.8]{nill}.

Now assume that $\delta_P>0$ and let's prove ($\ref{dueGM}$). 
For every vertex $v$ of $P$ we have $\langle v,u_i\rangle\geq
\delta_P$ if $i\not\in A(v)$;
similarly to what precedes,  we get
\begin{equation}
\label{pippo}
\begin{split}
0&=\sum_{i\in I}
m_i\langle v,u_i\rangle 
=-\sum_{i\in A(v)}m_i+\sum_{i\not\in A(v)} 
m_i\langle v,u_i\rangle \\&\geq -\sum_{i\in A(v)}m_i+\delta_P
\sum_{i\not\in A(v)} m_i
=(-\delta_P-1)\sum_{i\in A(v)}m_i+\delta_P M,
\end{split}
\end{equation}
namely $\frac{\delta_P}{\delta_P+1}M\leq\sum_{i\in A(v)}
m_i$. This gives
$$|V(P)|\frac{\delta_P}{\delta_P+1}M \leq \sum_{v\in V(P)}\sum_{i\in
  A(v)}m_i=nM\quad\text{and}\quad
|V(P)|\leq \frac{n (\delta_P+1)}{\delta_P}
=n+\frac{n}{\delta_P}.$$

Suppose that
$|V(P)|=n+\frac{n}{\delta_P}$. 
Again, we can choose a relation
as \eqref{sum}
involving all vertices of $P^*$.
Then, since we must have all equalities in \eqref{pippo}, we
conclude that $\langle
v,u\rangle\in\{-1,\delta_P\}$ for any
$v\in V(P)$ and $u\in V(P^*)$.

Set $r:=\frac{n}{\,\delta_P}$. Fix $u\in V(P^*)$ and call
$e_1,\dotsc,e_n$ the vertices of $F_{u}$,  
and $f_1,\dotsc,f_r$ the remaining vertices. 
Set $K:=\{1,\dotsc,n\}$ and $J:=\{1,\dotsc,r\}$.
For any $k\in K$, the face 
$\Conv(e_1,\dotsc,\check{e}_k,\dotsc,e_n)$
lies on exactly two facets, one of which is $F_{u}$. Hence there
exists a unique $\ph(k)\in J$ such that
$$F_k:=\Conv(f_{\ph(k)},e_1,\dotsc,\check{e}_k,\dotsc,e_n)$$ 
is a facet of $P$.
This defines a function $\ph\colon K\to J$.
 
Since $P$ is simplicial, $e_1,\dotsc,e_n$ is a basis of $N_{\Q}$; if
$e_1^*,\dotsc,e_n^*$ is the dual basis 
in $M_{\Q}$, we have
$u=-e_1^*-\cdots-e_n^*$.
Fix $k\in K$ and let 
$u_k$ be the vertex of $P^*$ such that $F_k=F_{u_k}$. 
We have 
$$\langle
e_i,u_k\rangle=-1    \text{ for all $i\in K\smallsetminus\{k\}\,$, and
}\langle e_k,u_k\rangle=\delta_P,$$ 
so  $u_k=u+(\delta_P+1)e_k^*$.

Now for any $j\in J$ we have
$$\langle f_j,e_k^*\rangle=\frac{1}{\delta_P+1}\left(\langle f_j,u_k\rangle
-\delta_P\right)=\begin{cases} -1 &\text{ if }\ph(k)=j,\\
                          0&\text{ otherwise.}\end{cases}$$
This means $f_j+\sum_{k\in\ph^{-1}(j)}e_k=0$.
Finally, we have $\delta_P=\langle f_j,u\rangle=|\ph^{-1}(j)|$, so as
$j$ varies in $J$, the  $\ph^{-1}(j)$'s give a partition of $K$ in $r$
subsets of cardinality $\delta_P$. Setting
$Q_j:=\Conv\{f_j,e_k\,|\,k\in\ph^{-1}(j)\}$, we see that
$Q_1,\dotsc,Q_r$ satisfy the properties claimed in ($\ref{dueGM}$).
\end{proof}
\begin{remark}
\label{relations}
Consider any polytope $Q\subset M_{\Q}$, of dimension $n$, containing
the origin in its interior. Let 
$u$ be a vertex of $Q$ and let
 $F$ be the minimal face of $Q$ such that $-u$ is contained in the
 cone over $F$ in $M_{\Q}$. Then writing $-u$ as a linear combination
 of the vertices of $F$, we get a relation as
 \eqref{sum} containing $u$. Summing enough 
relations of this type, one easily finds a 
relation  $\sum_{u\in V(Q)}m_uu=0$
with all $m_u$'s positive integers.

It is interesting to observe that when $Q$ is a reflexive polytope,
there is a special relation:
\begin{equation}
\label{folklore}
\sum_{u\in V(Q)}\Vol(F_u)\,u=0,\end{equation}
where $F_u$ is the facet of $Q^*$ corresponding to $u$, and 
$\Vol(F_u)$ is the lattice volume of $F_u$.
This follows from a theorem by
Minkowski, see \cite[page 332]{grunbaum} and \cite[Lemma 4.9]{nill2}.

If $P$ is a smooth reflexive polytope, 
then all facets of $P$ are standard simplices, so 
\eqref{folklore} yields that the sum of all vertices of $P^*$ is
zero. This remarkable fact can also be proven using the recent
results on the factorization of birational maps between smooth toric
varieties, and it was used for the proof of Theorem \ref{due} in a
previous version of this work.
\end{remark}

\smallskip

Let $P$ be a simplicial reflexive polytope. We denote by
 $\mathcal{N}_1(X_P)$  the $\Q$-vector space of 1-cycles  in $X_P$,
with rational
 coefficients, modulo numerical
 equivalence.  
 It is a well known fact in toric
geometry (see for instance \cite{torimori}) that there is an exact
 sequence
$$0\la \mathcal{N}_1(X_P)\la\Q^{V(P)}\la N_{\Q}\la 0,$$
so that $\mathcal{N}_1(X_P)$ is canonically identified
with the group of 
 linear rational relations among the vertices of $P$. 

Moreover, if a class $\gamma\in\mathcal{N}_1(X_P)$ 
corresponds to a relation
$$\sum_{v\in V(P)}m_vv=0,\qquad m_v\in\Q,$$
then the anticanonical degree of $\gamma$ is
$-K_{X_P}\cdot\gamma=\sum_{v\in V(P)}m_v$.
\begin{proof}[Proof of Lemma \ref{pseudoindex}]
To show that $\iota_{X_P}\leq\delta_P+1$, we 
exhibit a rational curve in $X_P$ whose anticanonical degree
is less or equal than $\delta_P+1$. 
Fix $v\in V(P)$ and $u\in V(P^*)$ 
such that $\langle v,u\rangle=\delta_P$. By Lemma \ref{zero}, $v$ is
adjacent to $F_u$. Let $e_1,\dotsc,e_n$ be the
vertices of $F_u$; up to reordering, 
we can assume that $\Conv(v,e_2,\dotsc,e_n)$ is a
facet of $P$. Since $P$ is simplicial, $e_1,\dotsc,e_n$ is a basis of
$N_{\Q}$; if $e_1^*,\dotsc,e_n^*$ is the dual
basis of $M_{\Q}$,
we have $u=-e_1^*-\cdots-e_n^*$.
Consider the relation
$$v+\sum_{i=1}^n a_i e_i=0,$$
and the corresponding class $\gamma\in\mathcal{N}_1(X_P)$. We have 
$$ -K_{X_P}\cdot \gamma=1+\sum_{i=1}^n a_i=1+\langle
v,u\rangle=1+\delta_P.$$
Now consider the invariant curve $C_0\subset X_P$ corresponding to the
cone over the face $\Conv(e_2,\dotsc,e_n)$. There exists
$b\in\Q$, $b\in(0,1]$, such that the numerical class of
$C_0$ is $b\gamma$ (see \cite[\S 2]{torimori}). 
Then $-K_{X_P}\cdot C_0=b(\delta_P+1)\leq\delta_P+1$. 

Assume now that $X_P$ is smooth, and
let $C$ be an invariant curve having minimal anticanonical
degree $\iota_{X_P}$. The numerical class of $C$ corresponds to a relation
$f_{0}+\sum_{i=1}^nb_if_i=0$, where $b_i\in\Z$ and
$\Conv(f_1,\dotsc,f_n)$ is a facet
$F$ of $P$ (see \cite[\S 2]{torimori}). 
The vertices $f_1,\dotsc,f_n$ are a basis of $N$; if
$f_1^*,\dotsc,f_n^*$ is the dual basis of $M$, then $F=F_{u_0}$ with 
$u_0=-f_1^*-\cdots-f_n^*$. So
$$ \iota_{X_P}=-K_{X_P}\cdot C=1+\sum_{i=1}^nb_i= 1+\langle
f_{0},u_0\rangle\geq 1+\delta_P\geq\iota_{X_P},$$ 
and $\iota_{X_P}=\delta_P+1$.
\end{proof}
\begin{proof}[Proof of Theorem \ref{uno}]
Part ($\ref{B}$) and the inequality in ($\ref{GM}$) are
straightforward consequences of Theorem \ref{due} and
Lemma \ref{pseudoindex}.

Assume that $\rho_{X_P}(\iota_{X_P}-1)=n$. Again using Lemma
\ref{pseudoindex} and Theorem \ref{due}, we get
$\delta_{P}\geq\iota_{X_P}-1>0$ and
$$|V(P)|\leq n+\frac{n}{\delta_P}\leq
n+\frac{n}{\iota_{X_P}-1}=n+\rho_{X_P}=|V(P)|,$$
so $\iota_{X_P}=\delta_P+1$, $|V(P)|=n+\frac{n}{\delta_P}$ and the
characterization in Theorem \ref{due} ($\ref{dueGM}$) holds.
This means that  $X_P$ is the quotient of $(\pr{\,\delta_P})^{\rho_{X_P}}$
by a finite subgroup $G$ of the big torus;
it is then enough to show that $X_P$ is smooth.

We keep the same notation as in the proof of Theorem \ref{due}. 
In order to show the smoothness of $X_P$, we have to show that
$e_1,\dotsc,e_n$ is a basis of $N$. Let $w\in N$ and write
$w=\sum_{i=1}^n\frac{t_i}{s_i}e_i$ with $t_i,s_i\in\Z$, $s_i\neq 0$,
and $t_i,s_i$ with no common factors for all $i$. 

Suppose that $t_1\neq 0$, we show that $s_1=1$ (for the other indices
the proof is analogous). Consider the relation
$$ f_{\ph(1)}+\sum_{k\in\ph^{-1}(\ph(1))}e_k=0,$$
and the corresponding class $\gamma\in\mathcal{N}_1(X_P)$. 
Now consider the invariant curve $C\subset X_P$
corresponding to the cone over the face $\Conv(e_2,\dotsc,e_n)$, and
recall that $\Conv(e_1,\dotsc,e_n)$ and
$\Conv(f_{\ph(1)},e_2,\dotsc,e_n)$ are faces of $P$.
Then there exists 
$b\in\Q$, $b\in(0,1]$, such that the numerical class of
$C$ is $b\gamma$ (see \cite[\S 2]{torimori}). So we have
$$\iota_{X_P}\leq -K_{X_P}\cdot C=b( -K_{X_P}\cdot \gamma)
=b(\delta_P+1)=b\iota_{X_P}
\leq\iota_{X_P},$$
which gives $b=1$. This is equivalent to saying that in the quotient lattice 
${N}/N\cap(\Q e_2\oplus\cdots\oplus\Q e_n)$, the image
$\overline{e}_1$ of $e_1$ is a generator. Now if $\overline{w}$ is the
image of $w$, we have $s_1\overline{w}=t_1\overline{e}_1$ with
$s_1,t_1$ non zero and with no common factors, so $s_1=1$.
\end{proof}

\medskip

\noindent {\bf Acknowledgments.} I am grateful to Benjamin
Nill for pointing out to me the 
relation \eqref{folklore}.

\bigskip

\begin{flushleft}
\textsc{C.\ Casagrande\\
Dipartimento di Matematica ``L.\ Tonelli'' \\
 Universit\`a di Pisa\\ 
Largo Bruno Pontecorvo, 5\\
56127 Pisa - Italy}\\
\textit{E-mail address}: casagrande@dm.unipi.it
\end{flushleft}

\footnotesize
\begin{thebibliography}{10}

\bibitem{occhettaGM}
Marco Andreatta, Elena Chierici, and Gianluca Occhetta.
\newblock Generalized {M}ukai conjecture for special {F}ano varieties.
\newblock {\em Central European Journal of Mathematics}, 2(2):272--293, 2004.

\bibitem{bat3}
Victor~V. Batyrev.
\newblock Toric {F}ano threefolds.
\newblock {\em Izvestiya Akademii Nauk SSSR Seriya Matematicheskaya},
  45(4):704--717, 1981 (in Russian).
\newblock English translation: Mathematics of the USSR Izvestiya, 19,
  13--25, 1982.

\bibitem{batdual}
Victor~V. Batyrev.
\newblock Dual polyhedra and mirror symmetry for {C}alabi-{Y}au hypersurfaces
  in toric varieties.
\newblock {\em Journal of Algebraic Geometry}, 3:493--535, 1994.

\bibitem{bat2}
Victor~V. Batyrev.
\newblock On the classification of toric {F}ano 4-folds.
\newblock {\em Journal of Mathematical Sciences (New York)}, 94:1021--1050,
  1999.

\bibitem{mukai}
Laurent Bonavero, Cinzia Casagrande, Olivier Debarre, and St{\'e}phane Druel.
\newblock Sur une conjecture de {M}ukai.
\newblock {\em Commentarii Mathematici Helvetici}, 78:601--626, 2003.

\bibitem{fano}
Cinzia Casagrande.
\newblock Toric {F}ano varieties and birational morphisms.
\newblock {\em International Mathematics Research Notices}, 27:1473--1505,
  2003.

\bibitem{cho}
Koji Cho, Yoichi Miyaoka, and Nick Shepherd-Barron.
\newblock Characterizations of projective space and applications to complex
  symplectic geometry.
\newblock In {\em Higher Dimensional Birational Geometry}, volume~35 of {\em
  Advanced Studies in Pure Mathematics}, pages 1--89. Mathematical Society of
  Japan, 2002.

\bibitem{debarreUT}
Olivier Debarre.
\newblock {\em Higher-Dimensional Algebraic Geometry}.
\newblock Universitext. Springer-Verlag, 2001.

\bibitem{debarre}
Olivier Debarre.
\newblock {F}ano varieties.
\newblock In {\em Higher Dimensional Varieties and Rational Points ({B}udapest,
  2001)}, volume~12 of {\em Bolyai Society Mathematical Studies}, pages
  93--132. Springer-Verlag, 2003.

\bibitem{ewald2}
G{\"u}nter Ewald.
\newblock {\em Combinatorial Convexity and Algebraic Geometry}, volume 168 of
  {\em Graduate Texts in Mathematics}.
\newblock Springer-Verlag, 1996.

\bibitem{grunbaum}
Branko Gr{\"u}nbaum.
\newblock {\em Convex Polytopes}, volume 221 of {\em Graduate Texts in
  Mathematics}.
\newblock Springer-Verlag, second edition, 2003.
\newblock First edition 1967.

\bibitem{nill2}
Benjamin Nill.
\newblock Complete toric varieties with reductive automorphism group.
\newblock {Preprint math.AG/0407491}, 2004.

\bibitem{nill}
Benjamin Nill.
\newblock Gorenstein toric {F}ano varieties.
\newblock {\em Manuscripta Mathematica}, 116(2):183--210, 2005.

\bibitem{occhetta}
Gianluca Occhetta.
\newblock A characterization of products of projective spaces.
\newblock {Preprint, available at the author's web page
  \verb1http://www.science.unitn.it/~occhetta/1}, 2003.

\bibitem{sato}
Hiroshi Sato.
\newblock Toward the classification of higher-dimensional toric {F}ano
  varieties.
\newblock {\em T{\^o}hoku Mathematical Journal}, 52:383--413, 2000.

\bibitem{VK}
V.~E. Voskresenski{\u{\i}} and Alexander Klyachko.
\newblock Toric {F}ano varieties and systems of roots.
\newblock {\em Izvestiya Akademii Nauk SSSR Seriya Matematicheskaya},
  48(2):237--263, 1984 (in Russian).
\newblock English translation: Mathematics of the USSR Izvestiya, 24,
  221--244, 1985.

\bibitem{wat}
Keiichi Watanabe and Masayuki Watanabe.
\newblock The classification of {F}ano 3-folds with torus embeddings.
\newblock {\em Tokyo Journal of Mathematics}, 5:37--48, 1982.

\bibitem{wisn2}
Jaros{\l}aw~A. Wi{\'s}niewski.
\newblock On a conjecture of {M}ukai.
\newblock {\em Manuscripta Mathematica}, 68:135--141, 1990.

\bibitem{torimori}
Jaros{\l}aw~A. Wi{\'s}niewski.
\newblock Toric {M}ori theory and {F}ano manifolds.
\newblock In {\em Geometry of Toric Varieties}, volume~6 of {\em S{\'e}minaires
  et Congr{\'e}s}, pages 249--272. Soci{\'e}t{\'e} Math{\'e}matique de France,
  2002.

\end{thebibliography}
\end{document}